\documentclass[10pt]{article}
\usepackage{amsmath,amssymb,amsthm}

\title{Splitting time for irrational triangle billiards}
\author{Dmitri Scheglov\\
University of Oklahoma}

\theoremstyle{plain}

\begin{document}

\maketitle
\setlength{\parindent}{0pt}

\begin{abstract}
\noindent
We give an upper estimate for the splitting time of a thin parallel beam of trajectories of irrational triangle billiard in terms of some number-theoretic function, depending on the angles of the billiard. 
\

In some cases we are able to explicitely estimate this function, thus giving an explicit upper estimate on the triangle billiard splitting time. 
\end{abstract}

\section{Introduction and main results}

In this paper we consider a billiard flow inside an irrational triangle. Here by irrational we mean that all the angles of the triangle are irrational multiples of $\pi$. As a well known fact in the theory of dynamical systems, there is a serious difference between billiard dynamics in rational and irrational polygons.
\

\

In the rational case one can apply the so-called Katok-Zemlyakov construction [5], which we will extensively use in this paper, and get a flat compact surface with singularities. The crucial observation here is that the abelian group , generated by rational billiard angles is finite, and so the phase space can be glued from the finite number of polygons. In this case the billiard flow is essentially a geodesic flow with singularities on this compact surface. 
\

\

At this point  we actually step into the realm of flat surfaces and Teichmuller theory comes into play. It serves as a very poverful working tool and allows to answer not only the basic questions, such as the complexity growtn, the growth of periodic orbits, directional ergodicity,[6],[7],[8]  but also more delicate as the deviation of ergodic averages in a given direction.[1]
\

\

This nice picture drastically changes as we consider billiards with irrational angles. First of all as the Katok-Zemliakov construction does not produce a compact surface, one can not apply Teichmuller theory any longer. And it turns out that most of the questions with established answers for the rational case are still open for irrational billiards. The main difficulty here is the absence of a tool, which would allow one to analyse the asymptotic behaviour of trajectories.
\

\

Nevertheless there are still a number of important results, describing different aspects of irrational billiard dynamics and we will mention some of them. 
\

R.Schwartz gave a rigorous computer-assisted proof that any triangle with angles less than 100 degrees has a periodic orbit. Let us also mention that for triangle with angles less than 90 degrees this fact is elementary.[9]
\

 Cipra, Hanson, Kolan proved that  for a right triangle a almost any orbit, perpendicular to a side, returns perpendicular and so is periodic.[2]
\

A.Katok proved that the complexity growth of orbits for any polygonal billiard is less than exponential. In fact it is conjectured that the growth is polynomial, since for rational billiards it is quadratic, but it is not clear how to get an explicit growth estimate, as Katok's proof is not explicit.[4]
\

Ya. Vorobets proved ergodicity of irrational billiards with a prescribed fast order of approximation by rational billiards by providing explicit estimates on close rational billiards. Unfortunately this idea does not work in general, because when one approximates irrational billiards by rational ones, then the genus of corresponding approximation surface grows and so ane has to have extrimely high speed of approximation to be able to provide explicit estimates. [10]
\

And finally the result most important for our purposes is proven by Gaplerin, Kruger, Troubeczkoy in [3]. In fact they have proven several important results in this paper, but the one, to which we pay a special attention is the following theorem:
\

\

\textbf{ Theorem ( Galperin, Kruger, Troubeczkoy).}
\

For any orbit of the polygonal billiard, not hitting any vertex, one of the two possibilities holds:
\

1) It is periodic.
\

2) Its closure contains one of the vertices.
\

\

Proof of this theorem is very elegant and uses abstract Furstenberg recurrence theorem. Since it uses an abstract fact, it is not constructive and in this paper in some sence we provide an explicit version of this theorem in case of irrational triangle, even though our formulation is slightly different.

\section{Definitions and the main theorem}
\

From here and further we will assume that the numbers $\alpha$, $\beta$, 1 are rationally independent, what corresponds to the generic billiard. In fact all analogous theorems would hold for the rationally dependent case and the proof would be even easier, but here we are interested in the  generic situation.
\

We consider the standard circle of radius 1 on the plane oriented counterclockwise with induced metric allowing obvious definition of the angluar difference $x-y$ between two points $x$ and $y$ on the circle and local linear structure near a fixed point.
\

From now and further using local linear structure of the circle we will naturally identify circle segments with linear intervals.
\

\

\textbf{Definition 1.} Consider a finite subset $S\subset\Delta$ of a segment $\Delta$ of the circle. $S$ is called a relative $\epsilon$-net if after linear "blowing up" of the segment $\Delta$ to the length 1, $S$ becomes an $\epsilon$-net in the standard sence.
\

\

\textbf{Definition 2.} Let $0<\alpha, \beta<\pi$. A finite sequence $x_{1},\ldots, x_{n}\in\mathbb{S}^{1}$ is called $\alpha\beta$-connected if for any $i: 1\leq n-1$ we have that either $x_{i+1}-x_{i}=\pm\alpha$ or $x_{i+1}-x_{i}=\pm\beta$ .
\

\

For any finite set of points $S\subset\mathbb{S}^{1}$  let $|S|$ denote its cardinality. We then have the following definition.
Fix a pair of numbers $\alpha$, $\beta$ as above. The $net$-function $ F_{\alpha \beta}$: (0,1)$\rightarrow Z_{+}$ is defined as follows.
\

\

\textbf{Definition 3.} $ F_{\alpha \beta}\left (\epsilon\right )=min \lbrace n\in\mathbb{Z}_{+}| $ $\forall  \alpha\beta $ -connected $S\subset\mathbb{S}^{1}$ $,$ $|S|>n$  $\exists $ $\overline S \subset S$, $\overline S$$-$ relative $\epsilon$-net  $\rbrace$
\

\

Informally speaking $ F_{\alpha \beta}\left (\epsilon\right )$ is a minimal cardinality of $\alpha\beta$ - connected sequence which guarantees that it contains a relative $\epsilon$-net.

\

We have to point out that it is an open question if the function  $ F_{\alpha \beta}$ exists for all $\alpha, \beta$. In this paper we will prove that it exists for some non-trivial pairs $\alpha, \beta$ and we conjecture the existence of a uniformly bounding function $F$ such that  $F_{\alpha\beta}\left (\epsilon\right )<F\left (\epsilon\right )$ for all small $\epsilon$.
\

\

\textbf{Definition 4.} Let $\alpha,\beta$ be rationally independent. Then $N_{\alpha\beta}\left ( k\right ) =$ min $\lbrace  \langle n\alpha+m\beta\rangle |$ for all $ |n|+|m|\leq k| \rbrace$, where $\langle x\rangle$ is a distance from $x$ to the closest integer.
\

\

Sometimes, when there is no ambiguity we will skip indices $\alpha\beta$ and use a notation $N\left (k\right ) $.
\

And we also give similar definition for one number.
\

\

\textbf{ Definition 5.} For an irrational number $\alpha$ let $N_{\alpha}\left (k\right )=$ min $\lbrace\langle n\alpha\rangle $ for all $|n|\leq k\rbrace$.
\

\

We will also use a notation $N\left ( k\right )$ when there is no ambiguity.
\

\

After all the needed definitions are made, we briefly remind to the reader a geometric unfolding construction, which from now and further we will call Katok-Zemlyakov construction.
\

 For a given polygon $P$ and any billiard trajectory $S$ inside $P$, instead of reflecting $S$ from the sides of $P$ and thus getting a piecewise linear billiard trajectory, we reflect polygon $P$ along the straight line and obtain a sequence of polygons. 
\

Because of the "mirror law" which governs the billiard trajectory, the segments of the straight line inside the different polygons along it are exactly the segments of the corresponding billiard trajectory.
\

The geometric shape resulting in the one time application of the Katok-Zemlyakov construction to a given triangle is called a KITE and is shown on the pic 1.
\begin{minipage}[c]{10mm}
\unitlength=1.00mm
\special{em:linewidth 0.4pt}
\linethickness{0.4pt}
\begin{picture}(23.67,44.67)
\put(13.00,44.67){\line(0,-1){41.67}}
\put(13.00,3.00){\line(-2,5){10.67}}
\put(2.34,30.00){\line(3,4){10.67}}
\put(13.00,44.34){\line(3,-4){10.67}}
\put(23.67,30.34){\line(-2,-5){11.00}}
\end{picture}

\end{minipage}

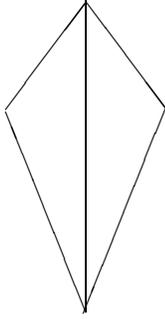
\textbf{Pic.1.} Kite, obtained as a result of one-time application of the Katok-Zemlyakov construction.
\

\

The next picture shows the Katok-Zemlyakov construction for a given billiard trajectory inside the kite.
\

\begin{minipage}[c]{10mm}
\unitlength=1.00mm
\special{em:linewidth 0.4pt}
\linethickness{0.4pt}
\begin{picture}(94.67,42.67)
\put(9.00,3.67){\line(-1,4){6.67}}
\put(2.33,29.67){\line(2,3){6.33}}
\put(8.66,39.34){\line(2,-3){6.67}}
\put(15.33,29.34){\line(-1,-4){6.00}}
\put(9.00,3.67){\line(1,1){17.33}}
\put(26.33,21.00){\line(0,1){13.00}}
\put(26.33,34.00){\line(-5,-2){11.33}}
\put(4.00,23.34){\vector(1,0){90.67}}
\put(26.33,34.00){\line(6,-1){11.00}}
\put(37.33,32.34){\line(1,-3){8.00}}
\put(45.33,8.34){\line(1,5){5.33}}
\put(45.33,8.34){\line(-3,2){19.00}}
\put(50.66,35.00){\line(-1,1){7.67}}
\put(43.33,42.67){\line(-3,-5){6.33}}
\end{picture}

\end{minipage}
\

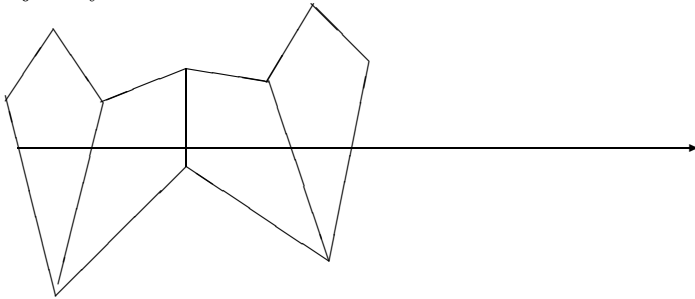
\textbf{Pic.2.} Katok-Zemlyakov construction for a kite.
\

\

We now fix a kite, obtained from a given triangle and we denote as $\alpha$ and $\beta$ the angles of a kite corresponding to the endpoints of the "reflecting" side of the triangle and assume that $\alpha<\beta$. On the Pic.1. $\alpha$ is a lower angle and $\beta$ is an upper angle.
\

Such a kite then has an oriented " main diagonal", connecting two vertices, corresponding to $\alpha$ and $\beta$ and we consider the sequence of angles $\theta_{i}$ between the main diagonal and straightened billiard trajectory as one applies Katok-Zemlyakov construction. We consider $\theta_{i}$ as a sequence of points on $\mathbb{S}^{1}$.
\

Simple and important observation which one immediately makes, looking at the pic.1. is the following:
\

$\theta_{i}$ is an $\alpha\beta$ - connected sequence!

\

Here we make one more useful geometric definition.
\

\textbf{Definition 6.} An $(\epsilon, T)$ - beam is a set of parallel segments, corresponding to the application of the Katok-Zemlyakov construction along some direction and from some base point, where $\epsilon$ is a width of the beam and $T$ is the length of the maximal parallel segment. 
\

The left interval on the kite side, transversal to the beam direction is called \underline {a base segment} or  \underline {base of a beam} and the right interval is called \underline{ an end segment} or \underline{end of the beam}.
\

\

Note, that by the definition beam does not have any kite vertices inside, as the Katok-Zemliakov is undefined on when the trajectory hits a vertex.
\

We will usually denote $(\epsilon, T)$ - beam as $B(\epsilon, T)$.
\

Below one may see a picture, which provides a geometric intuition behind the notion of $(\epsilon, T)$ - beam.
 \

\begin{minipage}[c]{10mm}
\unitlength=1.00mm
\special{em:linewidth 0.4pt}
\linethickness{0.4pt}
\begin{picture}(108.67,44.00)
\put(9.00,3.67){\line(-1,4){6.67}}
\put(2.33,29.67){\line(2,3){6.33}}
\put(8.66,39.34){\line(2,-3){6.67}}
\put(15.33,29.34){\line(-1,-4){6.00}}
\put(9.00,3.67){\line(1,1){17.33}}
\put(26.33,21.00){\line(0,1){13.00}}
\put(26.33,34.00){\line(-5,-2){11.33}}
\put(4.00,23.34){\vector(1,0){90.67}}
\put(26.33,34.00){\line(6,-1){11.00}}
\put(37.33,32.34){\line(1,-3){8.00}}
\put(45.33,8.34){\line(1,5){5.33}}
\put(45.33,8.34){\line(-3,2){19.00}}
\put(50.66,35.00){\line(-1,1){7.67}}
\put(43.33,42.67){\line(-3,-5){6.33}}
\put(91.67,12.00){\line(1,4){6.67}}
\put(3.00,27.00){\vector(1,0){92.67}}
\put(91.67,12.00){\line(4,5){17.00}}
\put(106.67,44.00){\line(-5,-3){8.67}}
\put(106.67,44.00){\line(1,-6){2.00}}
\put(108.67,32.67){\line(0,0){0.00}}
\end{picture}

\end{minipage}

\

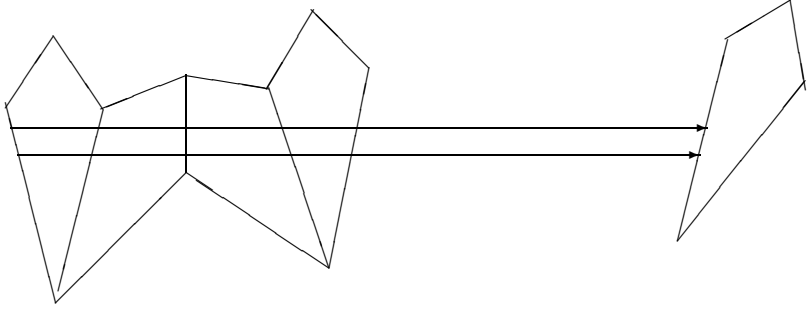
\textbf{Pic.3.} $(\epsilon, T)$ - beam of parallel trajectories.

\

\

We now prove several lemmas, which will be used later on.

\

\textbf{Lemma 1.} Let $\epsilon\in\left ( 0, 1/2\right )$ and $n$ be a positive integer. Let $\delta={\left (\frac{\epsilon}{100}\right )}^{n}$. Then for any relative $\delta$-net of an interval $\Delta$, colored in $n$ colors, there exists a monochromatic relative $\epsilon$-subnet of some subinterval $\overline{\Delta}\subseteq\Delta$.
\

\

\textbf{Proof.}. We will prove the claim by induction on $n$. As we work with relative nets, then without loss of generality we may assume $\Delta=\left [0, 1\right ]$.
\

If $n=1$ the claim is obvious. 
\

We now assume $n>1$ and divide $\Delta$ into $\frac{\epsilon}{2}$-segments by points $\frac{\epsilon}{2}, \epsilon, \frac{3}{2}\epsilon,\ldots, 1-p\epsilon$, where $p=\bigl [\frac{1}{\epsilon}\bigr ]$. If we can find a point of a  fixed color inside each segment, then we are done, because given points will form a monochromatic relative $\epsilon$-net on the interval $\left [0, 1-p\epsilon\right ]$.
\

  If not, there exists a segment $\sigma$ of width $\frac{\epsilon}{2}$ which has only the points of $n-1$ colors inside itself. Moreover there exists a segment $\overline\sigma\subset\sigma$ of width at least $\frac{\epsilon}{2}-2{\bigl (\frac{\epsilon}{100}\bigr )}^n$ with relative ${\left (\frac{\epsilon}{100}\right )}^{n}$-net colored in $n-1$ color. Blowing up $\overline\sigma$ to the length $1$ we get $\frac{{\left (\frac{\epsilon}{100}\right )}^{n}}{\frac{\epsilon}{2}-2{\bigl (\frac{\epsilon}{100}\bigr )}^n}<{\left (\frac{\epsilon}{100}\right )}^{n-1}$-net and the induction step is complete.
\

\

\textbf{Lemma 2.} Let $0<\gamma<\pi$ be an irrational number and $p, q\in\mathbb{Z}_{+}$. Let $\alpha=p\gamma, \beta=q\gamma$, and we also assume $\alpha, \beta<\pi$.
\

Let $S\subset\mathbb{S}^{1}$ be an $\alpha\beta$-connected sequence, $|S|>\frac{4\pi}{{\epsilon}^{2\left (p+q\right )}}\left (p+q\right ) 100^{2\left (p+q\right )}$ . 
\

Then $S$ as a subset of $\mathbb{S}^{1}$ contains a relative $\epsilon$ - net.
\

\textbf{Remark.$\left ( Important\right )$} The estimate in the lemma 2 is $\underline{ independent}$ on $\gamma$.

\

\textbf{Proof.} Thinking of $\mathbb{S}^{1}$ as of Abelian group and of $\alpha, \beta, \gamma $ as of its elements we note that $S\subset \gamma\mathbb{Z}\subset\mathbb{S}^{1}$. 
\

Here we assume that the first element $s_{1}\in S\subset\mathbb{S}^{1}$ corresponds to $0\in\mathbb{Z}$ and then $\gamma\mathbb{Z}$ as a subgroup of $\mathbb{S}^{1}$ is naturally isomorphic to $\mathbb{Z}$.

Since $\gamma$ is irrational then instead of $\alpha\beta$ - connected sequence $S$ we may also consider a $pq$ - connected sequence $S\subset\mathbb{Z}$, whose definition is completely analogous to that of $\alpha\beta$ - connected sequence. Note that we use the same symbol $S$ for both sequences as they completely define each other  and so it does not lead to ambiguity.
\

For any small enough $\delta>0$ there exists $n_{0}$ satisfying $p+q<n_{0}<2\pi\left (p+q\right )/\delta$ such that $|n_{0}\gamma$ mod $\mathbb{S}^{1}|<\delta$. It follows from the pigeonhole principle, applied to the sequence of points $0, \left (p+q\right )\gamma,\ldots, n\left (p+q\right )\gamma$ in $\mathbb{S}^{1}$. Then the sequence $0, n_{0}\gamma, \ldots, \bigl [\frac{1}{\delta}\bigr ]n_{0}\gamma$ forms a relative $\delta$ - net on the segment $\Delta=\bigl [0,\bigl [\frac{1}{\delta}\bigr ]n_{0}\gamma\bigr ]$.
\

We put $\delta= {\bigl (\frac{\epsilon}{100}\bigr )}^{p+q}$
\

Consider now $\alpha\beta$ connected sequence $S$ of cardinality $|S|>\frac{4\pi}{{\epsilon}^{2\left (p+q\right )}}\left (p+q\right ) 100^{2\left (p+q\right )}$. The corresponding sequence $S\subset\mathbb{Z}$ will then obviously have an element $s: |s|>\frac{2\pi}{{\epsilon}^{2\left (p+q\right )}}\left (p+q\right ) 100^{2\left (p+q\right )}$, which implies that either $s<-\frac{2\pi}{{\epsilon}^{2\left (p+q\right )}}\left (p+q\right ) 100^{2\left (p+q\right )} $ or $s> \frac{2\pi}{{\epsilon}^{2\left (p+q\right )}}\left (p+q\right ) 100^{2\left (p+q\right )}$.
\

Assume that $s>\frac{2\pi}{{\epsilon}^{2\left (p+q\right )}}\left (p+q\right ) 100^{2\left (p+q\right )}$. Second case is considered entirely analogously. Because $S\subset\mathbb{Z}$ is $pq$-connected it follows that it has at least one point $p_{k}$ in each of the intervals $\left [ kn_{0}, kn_{0}+p+q-1\right ]$, where $k$ runs through all the integers in the interval $\left [1,{\bigl (\frac{100}{\epsilon}\bigr )}^{p+q}\right ]$
\

To each of the points $kn_{0}$ we assign one of the $p+q$ "colors", namely a number $p_{k}-kn_{0}$.
\

Now we remind that points $\left [0, n_{0},\ldots, \left [{\bigl (\frac{100}{\epsilon}\bigr )}^{p+q}\right ]n_{0}\right ]$ form a relative $\delta$-net, considered as points on $\mathbb{S}^{1}$. By lemma 1 there exists an $\epsilon$-subnet $S_{0}$ of a fixed color.  By construction above the color is just a fixed number $t$ from the interval $\left [0, p+q-1\right ]$ which implies that the set $S_{0}+ t\gamma\subset\mathbb{S}^{1}$ is an $\epsilon$ - net from $S$, because rotations obviously preserve $\epsilon$ - nets.

\

\

\textbf{Remark.} As we will need the estimate from Lemma 2, we introduce the following useful notation:
\

$L\left (p,q,\epsilon\right )=\frac{4\pi}{{\epsilon}^{2\left (p+q\right )}}\left (p+q\right ) 100^{2\left (p+q\right )}$ 
\

\

\textbf{Theorem 1.} Let $\alpha, \beta$ be a pair of irrational numbers, such that $\alpha/\beta$ is irrational and moreover there exists a sequence of pairs $p_{n}, q_{n}\in\mathbb{Z}_{+}$ such that the following approximation inequality holds:
\

$|\alpha/\beta - p_{n}/q_{n}|<\frac{N_{\frac{\beta}{q_{n}}} \left (  L\left ( p_{n}, q_{n}, 1/n \right )                        \right )      }{100\beta{\left (100n\right )}^{L\left ( p_{n}, q_{n},1/n\right )} }$
\

For any small enough $\epsilon>0$ let us take $M=M\left (\epsilon\right )={\bigl ( \frac{100}{\epsilon}\bigr )}^{L\left ( p_{n},q_{n},\epsilon\right )}L\left ( p_{n}, q_{n},\epsilon\right )$, and $L=L\left (\epsilon\right )=L\left ( p_{n}, q_{n}, \epsilon\right )$, $R=R\left (\epsilon\right )= {\bigl ( \frac{100}{\epsilon}\bigr )}^{L\left ( p_{n},q_{n},\epsilon\right )}   $  where $n=\left [1/\epsilon\right ] + 1$. 
\

Then any $\alpha\beta$ - connected sequence $S$ such that $|S|>M\left (\epsilon\right )$ contains a relative $\epsilon$ - net.
\

\

\textbf{Remark.} This theorem implies that for a pair of numbers $\alpha, \beta$ as in the theorem assumption, we have $F_{\alpha\beta}\left (\epsilon\right )\leq M\left (\epsilon\right )$ and in particular it means that $F_{\alpha\beta}$ exists in this case.

\

\textbf{Proof}. We first introduce some notations. We fix small enough $\epsilon>0$ and let $n, p_{n}, q_{n}$ be defined as in the formulation of the theorem. Let $\gamma=\beta/q_{n}$ and $\overline\alpha=p_{n}\gamma$ and of course $\beta=q_{n}\gamma$. Let also $S_{T}\subset\mathbb{S}^{1}$ be defined as $S_{T}=\lbrace k\gamma, k\in\left [-T, T\right ]\rbrace$. 
\

Let $\mu=\frac{N_{\frac{\beta}{q_{n}}} \left (  L\left ( p_{n}, q_{n}, 1/n \right )                        \right )      }{100\beta }$ and let and let $U_{T}$ be a $\mu$ - neighbourhood of $S_{T}$ for any $T\in\mathbb{Z}_{+}$.
\

For any $\alpha\beta$ - connected sequence $S\subset U_{L}$ we have a uniquely defined $\overline\alpha\beta$ - connected sequence $\overline S$.
\

Let us now take a closer look of on the $\alpha\beta$ - connected sequence $S$ assuming for a moment that $S\subset U_{L-1}$.
\

Let $s_{i}\in S$ and by our assumption there is a uniquely defined close point $s\in S_{T}$ such that the distance dist$\left (s_{i}, s\right )=$ dist$\left (s_{i}, S\right )$.
\

Now we remind that we have orientation on the circle and consider the walue $w_{i}=s_{i}-s$ with respect to the orientation. As $S$ is $\alpha\beta$ - connected, then $s_{i+1}=s_{i}\pm\alpha$ or $s_{i+1}=s_{i}\pm\beta$, which means that either $w_{i+1}=w_{i}$ or $w_{i+1}=w_{i}\pm\delta$ where $\delta=\alpha-\overline\alpha$.
\

Informally speaking this means that "from the point of view of the set $U_{M}$" we either do not move or make a shift on $\pm\delta$. 
\

More formally the sequence $w_{i}$ is $\delta$ - connected, where the definition of $\delta$ - connectedness is completely analogous to that of $\alpha\beta$ - connectedness.
\

We now compare three sequences: $\alpha\beta$ - sequence $S=\left (S_{i}\right )$, $\overline\alpha\beta$ - sequence $\overline S=\left (\overline S_{i}\right )$ and $\delta$ - sequence $w=\left (w_{i}\right )$.
\

Assume that $|S|>M\left (\epsilon\right )$ and $S$ does not contain a relative $\epsilon$-net. Then two situations are possible.
\

\

\textbf{Case 1.} $S\subset U_{L}$. Then, since $|S_{M}|<L$ and $|S|>M$ by pigeonhole principle $|w|>R$. We now consider "colors", corresponding to the points of $\overline S_{L}$ and we see that the sequence $w$ is naturally colored in at most $L$ colors.  Since $w$ is $\delta$ - connected and $|w|>R$, then $w$ as a set forms a relative $\mu$ - net, colored in $L$ color.
\

By Lemma 1 there exists a relative subnet of one color. But by the construction of colors this precisely means that $S$ has a relative $\epsilon$ - net, located near one of the points of the set $\overline S_{N}$ and so in this case we obtain a contradiction.
\

\

\textbf{Case 2.} $S\nsubseteq U_{L}$. In this case consider the largest index $n$ such that the sequences $S_{n}$ of the first $n-1$ terms $s_{1},\ldots, s_{n-1}$ belong to $U_{L}$. We also consider a corresponding terms of a sequence $\overline S$ and let $\overline S_{n}$ be a sequence $\overline s_{1},\ldots, \overline s_{n}$ and and let $\overline W_{n}$ be a sequence $w_{1},\ldots, w_{n}$  of corresponding terms of a sequence $w$  .  As $s_{n}\notin U_{L}$ then either $|\overline S_{n}|>L$ or $|\overline S_{n}|\leq L$ ,     $|M|>N$. 
\

In the first case $S$ contains a relative $\epsilon$-net  because $\overline S_{n}$ is $\overline \alpha\beta$ - connected the claim follows from Lemma.

\

In the second case by pigeonhole principle we have $|\overline W_{n}|>R$ and each point from $\overline W_{n}$ is naturally colored in not more then $L$ colors and the claim follows from the "color"argument used before and Lemma.
\

\

The following technical lemma will be used in our proof of the theorem 2.

\

\textbf{Lemma 3.} Let $B(\epsilon, T)$ be a parallel beam such that :
\

$1 )$ The beam base intersects with a beam end as kite segments.
\

$2 )$ The base and end angular kite positions are the same.
\

Then there is a segment inside the beam, corresponding to the periodic orbit.
\

\begin{minipage}[c]{10mm}
\unitlength=1.00mm
\special{em:linewidth 0.4pt}
\linethickness{0.4pt}
\begin{picture}(115.67,91.33)
\put(3.33,52.33){\line(2,3){15.00}}
\put(18.33,75.00){\line(2,-3){15.67}}
\put(85.00,68.67){\line(2,3){15.00}}
\put(100.00,91.33){\line(2,-3){15.67}}
\put(10.33,30.67){\line(1,0){87.33}}
\put(3.66,50.67){\line(1,0){87.67}}
\put(10.00,32.00){\line(5,1){82.00}}
\put(3.33,52.00){\line(1,-3){15.67}}
\put(19.00,5.33){\line(1,3){15.33}}
\put(19.00,5.33){\line(3,2){36.67}}
\put(55.66,29.33){\line(1,6){4.67}}
\put(85.00,68.33){\line(1,-3){15.33}}
\put(100.33,22.67){\line(1,3){15.00}}
\put(34.00,51.00){\line(4,1){26.33}}
\end{picture}

\end{minipage}
\

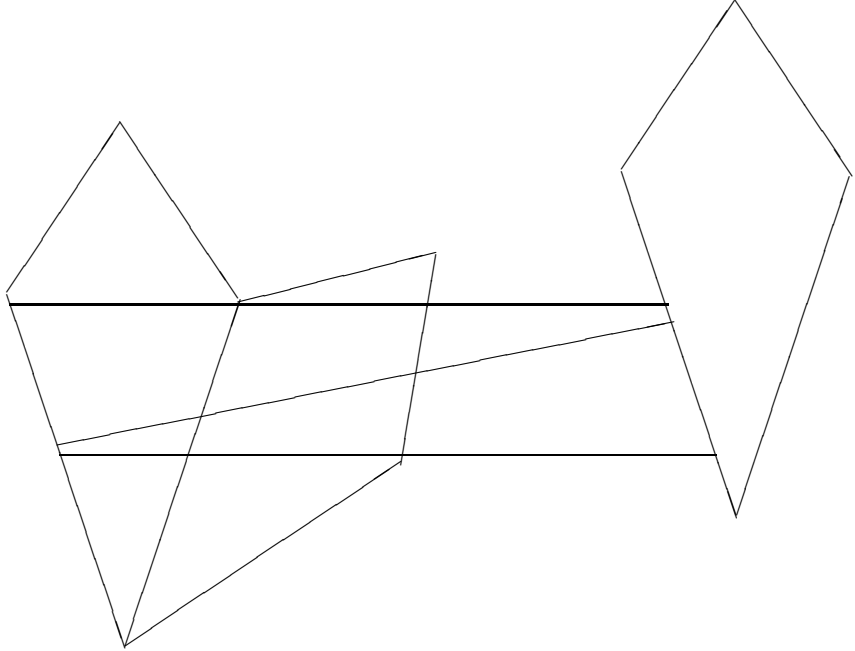
\textbf{Pic.5.} Periodic trajectory inside a beam with intersecting base and end.

\

\textbf{Proof.} As the base and end of the beam intersect, then there is a straight segment connecting them. And as the angular positions are the same, this straight segment corresponds to the periodic orbit. See the picture above for the geometric intuition behind the proof.

\

\textbf{Lemma 4.} Let $B(\epsilon, T)$ be a beam, where $\epsilon$ is small enough. Then the number of kites it intersects during the Katok-Zemlyakov construction is less then $CT/\epsilon$, where $C$ is a constant depending on kite only.
\

\

\textbf{Proof.} Consider any intersection of the beam with a given kite. Since the width of the beam is $\epsilon$ then for $\epsilon$ small enough there will be at least on point of the beam on the distance at least $\epsilon$ from all the kite vertices. Then by elementary trigonometry there exists a beam segment inside a kite of the length at least $\epsilon/C$, where $C>0$ is large enough and depend only on the kite.
\

 Now counting the total length of the beam segments for all kites intersecting the beam gives the desired inequality.

\

We now have all the tools to prove the main theorem.
\

\

\textbf{Theorem 2.}  Let $\alpha$, $\beta$ be a pair of rationally independent numbers such that $ F_{\alpha \beta}$($\epsilon$) is a correctly defined function for positive $\epsilon$. And let  $ K_{\alpha \beta}$ be the kite of diameter 1 with angles  $\alpha$, $\beta$. Let B($\epsilon$, $M$) be a parallel beam. 
\

Let us also introduce the following notations: $P_{\alpha\beta} (\epsilon) = F_{\alpha\beta} \bigl ( {\bigl ( \frac{\epsilon}{1600} \bigr )}^{[ 16/\epsilon ]+1}            \bigr )$
\

$Q_{\alpha\beta}(\epsilon)=\bigl (\bigl [\frac{16}{\epsilon}\bigr ] + 1\bigr )P_{\alpha\beta}(\epsilon)$
\

Then either $B$ as a set contains a periodic trajectory or $M\leq Q_{\alpha\beta}(\epsilon)+\frac{C}{\epsilon N_{\alpha\beta} (Q_{\alpha\beta}(\epsilon))}$,
where $C$ is a constant, depending only on the kite.

\

\textbf{Proof.} We first fix $\epsilon>0$ and divide the perimeter of our kite into $\bigl [\frac{16}{\epsilon}\bigr ]+1$ segments , so that each segment will be of the length less than $\frac{\epsilon}{4}$. We assign a specific "color" to each segment, which is just a positive integer from $1$ to $\bigl [\frac{16}{\epsilon}\bigr ]+1$.
\

We then take a middle point $x$ in the beam base and consider its trajectory $x_{1}$,$\ldots$,$x_{n}$ under the discrete billiard map. And we also consider a corresponding sequence of based vectors $v_{1},\ldots,v_{n}$ in the phase spase of the discrete billiard as shown on the pic.6.
\

\begin{minipage}[c]{10mm}
\unitlength=1.00mm
\special{em:linewidth 0.4pt}
\linethickness{0.4pt}
\begin{picture}(26.00,55.07)
\put(14.00,4.67){\line(-1,3){12.00}}
\put(2.00,40.67){\line(5,6){12.00}}
\put(14.00,4.67){\line(1,3){12.00}}
\put(26.00,40.67){\line(-5,6){12.00}}
\put(9.66,18.00){\vector(1,2){4.67}}
\put(5.33,30.00){\vector(2,1){8.33}}
\put(7.00,46.34){\vector(1,0){9.00}}
\put(23.33,32.67){\vector(-1,1){7.00}}
\put(3.00,38.00){\vector(1,0){8.67}}
\put(19.00,20.00){\vector(0,1){10.00}}
\end{picture}

\end{minipage}

\

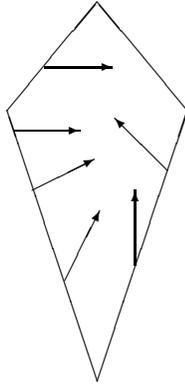
\textbf{Pic.5.} Kite and a sequence of vectors for a discrete billiard map.
\

\

  As we assigned a color to each segment, we may say that the two sequences above become " colored". In case if the trajectory hits on of the dividing points, we assign any color from one of the two segments.
\

\

As the unit tangent bundle to the kite is naturally trivial, we have a sequence of points on the unit circle $S_{n}=s_{1},\ldots,s_{n}$ corresponding to the sequence $v_{1},\ldots,v_{n}$. We remind that the sequence $S_{n}$ is $\alpha\beta$-connected.
\

Now assume that $n>Q_{\alpha\beta}(\epsilon)$. We have two possibilities on the cardinality of $S_{n}$.
\

\

\textbf{1 case.} $|S_{n}|< P_{\alpha\beta}(\epsilon)$ 
\

\

Then by the pigenhole principle we have at least $\bigl (\big [\frac{16}{\epsilon}\bigr ] + 2\bigr )$ points of $M$ equal to a given fixed point $s\in\mathbb S^{1}$.
Since our points are colored in $\bigl (\big [\frac{16}{\epsilon}\bigr ] + 1\bigr )$ colors, applying pigeonhole principle again we see that there are at least two points $s_{i},s_{j}=s\in M$ of the same color.
\

Since the color is the same it follows that the corresponding points $x_{i}, x_{j}$ both lie in the same segment of width at most $\epsilon/4$. Since $s_{i}=s$ and $s_{j}=s$ it implies that the corresponding vectors $v_{i}$ and $v_{j}$ are parallel. But as the width of the beam is $\epsilon$ it means that we have all the assumptions of the Lemma 1, and so there is a periodic orbit inside.
\

\begin{minipage}[c]{10mm}
\unitlength=1.00mm
\special{em:linewidth 0.4pt}
\linethickness{0.4pt}
\begin{picture}(125.34,86.60)
\put(9.67,2.34){\line(-1,3){8.33}}
\put(1.33,27.34){\line(5,6){7.67}}
\put(9.00,36.67){\line(5,-6){7.33}}
\put(16.33,27.67){\line(-1,-4){6.33}}
\put(6.00,14.67){\line(6,1){119.00}}
\put(5.67,14.67){\line(4,1){119.33}}
\put(5.67,14.67){\line(3,1){119.67}}
\put(5.33,15.34){\line(5,2){119.67}}
\put(5.67,14.67){\line(2,1){119.67}}
\put(5.67,15.00){\line(5,3){119.33}}
\put(102.33,36.67){\line(1,4){7.00}}
\put(124.67,34.34){\line(0,1){52.00}}
\put(102.33,36.67){\line(4,3){21.00}}
\put(123.33,52.33){\line(-1,5){2.33}}
\put(109.33,64.33){\line(1,0){12.00}}
\end{picture}

\end{minipage}
\

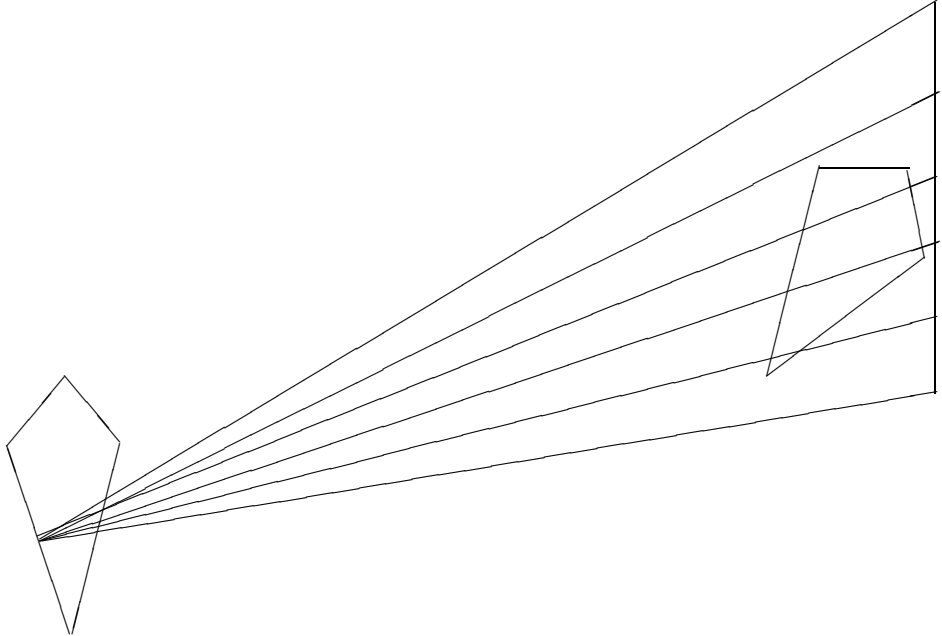
\textbf{Pic.6.} The " long triangle" formed by $\epsilon$ - net and a kite, "caught" in the $\epsilon$ - net in the case 2.

\

\textbf{2 case.}  $|S_{n}|\geq P_{\alpha\beta}(\epsilon)$ 
\

Then by definition of the function $F_{\alpha\beta}$ there exists a relative  $\epsilon/16$-net of the same color, which can be reformulated as follows. There exists a finite set of points $y_{1},\ldots, y_{k}$ from the trajectory of $x$, lying in the same segment $\sigma$, such that the corresponding vectors $w_{1},\ldots,w_{k}$ form a relative $\epsilon/16$-net on the circle $\mathbb S^{1}$.
\

\

Since all the vectors above belong to the same $x$-trajectory, we may start unfolding from each of the points  $y_{1},\ldots, y_{k}$ in the corresponding directions $w_{1},\ldots,w_{k}$. Since by assumption we do not hit a vertex during the time, the unfolding procedure goes exactly the same way for all the vectors $w_{1},\ldots,w_{k}$ as is shown on the pic.6. 
\

But since the diameter of the kite $ diam\bigl (K\bigr)<1$ we have an obvious contradiction, as the strips of width $\epsilon$ completely cover the corresponding " long triangle", with the right base equal to 1 and so one of the strips must hit one of the kite vertices. As the " width of each strip is less then $\epsilon$ we get a contradiction.
\

\

Let us now estimate the maximal length of the " long rays" on Pic.6.  As the total number of iterations we have made is $Q_{\alpha\beta}(\epsilon)$ then the angles between the rays are at least $N_{\alpha\beta}(Q_{\alpha\beta}(\epsilon))$. As the base of the " long triangle" is 1, and for small enough $x$ , we have $sin(x)\approx x$, then the maximal ray length $L\leq\frac{C}{N_{\alpha\beta}(Q_{\alpha\beta}(\epsilon))}$.
\

Then due to the Lemma 4 the number of kites intersecting the maximal length ray is at most $ \frac{C} {\epsilon N_{\alpha\beta} ( Q_{\alpha\beta}(\epsilon))}$
\

Adding the last estimate with the upper estimate on the time, needed to get a "long triangle" $Q_{\alpha\beta}(\epsilon)$ we complete the proof.

\

\

\textbf {3. Applications }.
\

In the upcoming paper in preparation we apply the estimate from the theorem 2 and partially the tools used in the proof of the theorem 2 to obtain lower bounds on the directional complexity in all non-periodic and almost all directions for kites with angles $\alpha, \beta$ for which the function $F_{\alpha\beta}(\epsilon)$ is correctly defined.
\

\

\textbf{Open problems.}
\

\

1) As the current paper has shown, the function $F_{\alpha\beta}(\epsilon)$ allows to obtain essential information about the billiard and in particular , about its directional complexity, which will be shown in the upcoming paper. Which allows us to formulate the following important problem:
\

\

\textbf{Problem 1.} Prove that $F_{\alpha\beta}(\epsilon)$ exists for any pair of irrational $\alpha\beta$.
\

\

The second problem is a stronger conjecture based on some heuristic observations about the behaviour of $F_{\alpha\beta}(\epsilon)$
\

\

\textbf{Problem 2.} Prove the following universal estimate for all irrational $\alpha,\beta$:

\

$F_{\alpha\beta}(\epsilon)<F(\epsilon)$
\

\

where $F(\epsilon)$ is a function \underline{independent} on $\alpha, \beta$.

\end{document}